\documentclass[12pt]{amsart}
\usepackage{amsmath, amssymb}
\numberwithin{equation}{section}

\theoremstyle{plain}
\newtheorem{theorem}{Theorem}[section]

\newtheorem{lemma}[theorem]{Lemma}

\theoremstyle{definition}
\newtheorem{remark}[theorem]{Remark}
\newtheorem{definition}[theorem]{Definition}

\def \N {\mathbb{N}}
\def \R {\mathbb{R}}

\def \E {\mathbb{E}}
\def \P {\mathbb{P} \, }

\def \a {\alpha}
\def \b {\beta}

\def \e {\varepsilon}
\def \d {\delta}
\def \D {\Delta}
\def \f {\varphi}

\def \s {\sigma}
\def \t {\tau}
\def \om {\omega}
\def \Om {\Omega}

\def \etc {,\ldots,}

\def \vol {{\rm vol}}

\def \supp {{\rm supp}}

\newcommand{\sumn}{\sum_{j=1}^n}
\newcommand{\pr}[2]{\langle {#1} , {#2} \rangle}
\newcommand{\norm}[1]{\left \| #1 \right \|}
\newcommand{\cref}[1]{c_{\text{\rm \ref{#1}}}}
\newcommand{\Cref}[1]{C_{\text{\rm \ref{#1}}}}

\begin{document}
\title {Invertibility of random matrices:
        norm of the inverse}
\author {Mark Rudelson}
\address{Department of Mathematics \\
   University of Missouri \\
   Columbia, MO 65211. }
\email{rudelson@math.missouri.edu}
\urladdr{http://math.missouri.edu/$\sim$rudelson}
\thanks{Research was supported in part by  NSF grant DMS-024380.}

\begin{abstract}
  Let $A$ be an $n \times n$ matrix, whose entries are independent
  copies of a centered random variable satisfying the subgaussian tail
  estimate. We prove that the operator norm of $A^{-1}$ does not exceed
  $Cn^{3/2}$ with probability close to $1$.
\end{abstract}

\maketitle

\section{Introduction.}   \label{s: introduction}

Let $A$ be an $n \times n$ matrix, whose entries are independent
identically distributed random variables. The spectral properties
of such matrices, in particular invertibility, have been
extensively studied (see, e.g. the survey \cite{DS}). While $A$ is
almost surely invertible whenever its entries are absolutely
continuous, the case of discrete entries is highly non-trivial.
Even in the case, when the entries of $A$ are independent random
variables taking values $\pm 1$ with probability $1/2$, the
precise order of probability that $A$ is degenerate is unknown.
Koml\'{o}s \cite{K1, K2} proved that this probability is $o(1)$ as
$n \to \infty$. This result was improved by Kahn, Koml\'{o}s ans
Szemer\'{e}di \cite{KKS}, who showed that this probability is
bounded above by $\theta^n$ for some absolute constant $\theta<1$.
The value of $\theta$ has been recently improved in a series of
papers by Tao and Vu \cite{TV1, TV2} to $\theta=3/4+o(1)$ (the
conjectured value is $\theta=1/2+o(1)$).

 However, these papers do
not address the quantitative characterization of invertibility,
namely the norm of the inverse matrix, considered as an operator
from $\R^n$ to $\R^n$. Random matrices are one of the standard
tools in geometric functional analysis. They are used, in
particular, to estimate the Banach--Mazur distance between
finite-dimensional Banach spaces and to construct sections of
convex bodies possessing certain properties. In all these
questions the distortion $\norm{A} \cdot \norm{A^{-1}}$ plays the
crucial role. Since the norm of $A$ is usually highly
concentrated, the distortion is determined by the norm of
$A^{-1}$. The estimate of the norm of $A^{-1}$ is known only in
the case when $A$ is a matrix with independent $N(0,1)$ Gaussian
entries. In this case Szarek \cite{Sz2} proved that $\norm{A^{-1}}
\le c\sqrt{n}$ with probability close to $1$ (see also \cite{Sz1}
where the spectral properties of a Gaussian matrix are applied to
an important question from  geometry of Banach spaces). For other
random matrices,including a random $\pm 1$ matrix, even a
polynomial bound was unknown. Proving such polynomial estimate is
the main aim of this paper.

More results are known about rectangular random matrices. Let
$\Gamma$ be an $N \times n$ matrix, whose entries are independent
random variables. If $N>n$, then such matrix can be considered as
a linear operator $\Gamma:\R^n \to Y$, where $Y=\Gamma \R^n$. If
we consider a family $\Gamma_n$ of such matrices with $n/N = \a$
for a fixed constant $\a>1$, then the norms of
$(\Gamma_n|_Y)^{-1}$ converge a.s. to $(1-\sqrt{\a})^{-1}
n^{-1/2}$, provided that the fourth moments of the entries are
uniformly bounded \cite{BY}. The random matrices for which
$n/N=1-o(1)$ are considered in \cite{LPRT}. If the entries of such
matrix satisfy certain moment conditions and $n/N > 1- c/\log n$,
then $\norm{(\Gamma|_Y)^{-1}} \le C(n/N) \cdot  n^{-1/2}$ with
probability exponentially close to $1$.

The proof of the last result is based on the $\e$-net argument. To
describe it we have to introduce some notation. For $p \ge 1$ let
$B_p^n$ denote the unit ball of the Banach space $\ell_p^n$. Let
$E \subset \R^n$ and let $B \subset \R^n$ be a convex symmetric
body. Let $\e>0$. We say that a set $F \subset \R^n$ is an
$\e$-net for $E$ with respect to $B$ if
\[
  E \subset \bigcup_{x \in F} (x+ \e B).
\]
The smallest cardinality of an  $\e$-net will be denoted by $N(E,
B, \e)$. For a point $x \in \R^n$, $\norm{x}$ stands for the
standard Euclidean norm, and for a linear operator $T: \R^n \to
\R^m$, $\norm{T}$ denotes the operator norm of $T: \ell_2^n \to
\ell_2^m$.

Let $E \subset S^{n-1}$ be a set such that for any fixed $x \in E$
there is a good bound for the probability that $\norm{\Gamma x}$
is small. We shall call such bound the small ball probability
estimate. If $N(E, B_2^n, \e)$ is small, this bound implies that
with high probability $\norm{\Gamma x}$ is large for all $x$ from
an $\e$-net for  $E$. Then the approximation is used to derive
that in this case $\norm{\Gamma x}$ is large for all $x \in E$.
Finally, the sphere $S^{n-1}$ is partitioned in two sets for which
the above method works. This argument is possible because the
small ball probability is controlled by a function of $N$, while
the size of an $\e$-net depends on $n<N$.

The case of a square random matrix is more delicate. Indeed, in
this case the small ball probability estimate is too weak to
produce a non-trivial estimate for the probability that
$\norm{\Gamma x}$ is large for all points of an $\e$-net. To
overcome this difficulty, we use the $\e$-net argument for one
part of the sphere and work with conditional probability on the
other part. Also, we will need more elaborate small ball
probability estimates, than those employed in \cite{LPRT}. To
obtain such estimates we use the method of Hal\'{a}sz, which lies
in the foundation of the arguments of \cite{KKS}, \cite{TV1},
\cite{TV2}.

Let $\P(\Om)$ denote the probability of the event $\Om$, and let
$\E \xi$ denote the expectation of the random variable $\xi$. A
random variable $\b$ is called subgaussian if for any $t>0$
\begin{equation}  \label{subgaussian}
  \P(|\b|>t) \le C \exp (-ct^2).
\end{equation}
The class of subgaussian variables includes many natural types of
random variables, in particular, normal and bounded ones.  It is
well-known that the tail decay condition \eqref{subgaussian} is
equivalent to the moment condition $\big (\E|\b|^p \big )^{1/p}
\le C' \sqrt{p}$ for all $p \ge 1$.

The letters $c, C, C'$ etc.  denote unimportant absolute
constants, whose value may change from line to line. Besides these
constants, the paper contains many absolute constants which are
used throughout the proof. For reader's convenience we use a
standard notation for such important absolute constants. Namely,
if a constant appears in the formulation of Lemma or Theorem x.y,
we denote it $C_{\text{x.y}}$ or $c_{\text{x.y}}$.

The main result of this paper is the following polynomial bound
for the norm of $A^{-1}$.

\begin{theorem}  \label{t: main}
  Let $\b$ be a centered subgaussian random variable of variance
  1.
  Let $A$ be an $n \times n$ matrix whose entries are independent
  copies of $\b$. Then for
  any $\e>\cref{t: main}/\sqrt{n}$
  \[
    \P \left ( \exists \, x \in \R^n \mid
               \norm{Ax}< \frac{\e}{\Cref{t: main} \cdot n^{3/2}}
       \right )
    < \e
  \]
  if $n$ is large enough.
\end{theorem}
More precisely, we prove that the probability above is bounded by
$\e/2 + 4 \exp( - c n )$ for all $n \in \N$.

The inequality of Theorem \ref{t: main} means that $\norm{A^{-1}}
\le \Cref{t: main} \cdot n^{3/2}/\e$ with probability greater than
$1-\e$. Equivalently, the smallest singular number of $A$ is at
least $\e / (\Cref{t: main} \cdot n^{3/2})$.

An important feature of Theorem \ref{t: main} is its universality.
Namely, the probability estimate holds for all subgaussian random
variables, regardless of their nature. Moreover, the only place,
where we use the assumption that $\b$ is subgaussian, is Lemma
\ref{l: norm of A} below.

\section{Preliminary results.}   \label{s: preliminary}

Assume that $l$ balls are randomly placed  in $k$ urns. Let $V \in
\{1 \etc k\}^l$ be a random vector whose $i$-th coordinate is the
number of balls contained in the $i$-th urn. The distribution of
$V$, called random allocation, has been extensively studied, and
many deep results are available (see \cite{KSC}). We need only a
simple combinatorial lemma.
\begin{lemma}    \label{l: concentration}
  Let $k \le l$
  and let $X(1) \etc X(l)$ be i.i.d. random variables
  uniformly distributed on the set $\{1 \etc k\}$.
  Let $\eta<1/2$.
  Then with probability greater than $1-\eta^l$ there exists a
  set $J \subset \{1 \etc l\}$ containing at least $l/2$ elements
  such that
  \begin{equation}  \label{profile for I}
     \sum_{i=1}^k |\{j \in J \mid X(j)=i \}|^2
     \le C(\eta) \frac{l^2}{k}.
  \end{equation}
\end{lemma}
\begin{remark}
  The proof yields $C(\eta)= \eta^{-16}$. This estimate is by no
  means exact.
\end{remark}
\begin{proof}
Let $X=(X(1) \etc X(l))$. For $i=1 \etc k$ denote
\[
  P_i(X)=  |\{j \mid X(j)=i\}|.
\]
 Let $2<\a<k/2$ be a number to be chosen later. Denote
\[
  I(X)=\{i \mid P_i(X) \ge \a \, \frac{l}{k} \}.
\]
For any $X$ we have $\sum_{i=1}^k P_i(X)=l$, so $|I(X)| \le k/\a$.
Set $J(X)=\{ j \mid X(j) \in I(X) \}$.  Assume that $|J(X)| \le
l/2$. Then for the set $J'(X)=\{1 \etc l\} \setminus J(X)$ we have
$|J'(X)| \ge l/2$ and
\[
  \sum_{i=1}^k |\{j \in J'(X) \mid X(j)=i\}|^2
  =\sum_{i \notin I(X)} P_i^2(X)
  \le k \cdot \left( \a \, \frac{l}{k} \right )^2
  = \frac{\a^2 l^2}{k}.
\]
Now we have to estimate the probability that $|\{J(X)\}| \ge l/2$.
To this end we estimate the probability that $J(X)=J$ and $I(X)=I$
for fixed subsets $J \subset \{1 \etc l\}$ and $I \subset \{1 \etc
k\}$ and sum over all relevant choices of $J$ and $I$. We have
\begin{align*}
  \P (|J(X)| \ge l/2)
  &\le \sum_{|J|\ge l/2} \ \sum_{|I| \le k/\a} \P (J(X)=J, \ I(X)=I)
  \\
  &\le \sum_{|J|\ge l/2} \ \sum_{|I| \le k/\a}
                        \P (X(j) \in I \text{ for all } j \in J) \\
  &\le 2^l (k/\a) \cdot \binom{k}{k/\a} \cdot (1/\a)^{l/2} \\
  &\le k \cdot (e \a)^{k/\a} \cdot (4/\a)^{l/2},
\end{align*}
since the random variables $X(1) \etc X(l)$ are independent. If $k
\le l$ and $\a>100$, the last expression does not exceed
$\a^{-l/8}$. To complete the proof, set $\a=\eta^{-8}$. If
$\eta>(2/k)^{1/8}$, then the assumption $\a<k/2$ is satisfied.
Otherwise, we can set $C(\eta)=\a^2 >(k/2)^2$, for which the
inequality \eqref{profile for I} becomes trivial.

\end{proof}

The following result is a standard large deviation estimate (see
e.g. \cite{DS} or \cite{LPRT}, where a more general result is
proved).
\begin{lemma} \label{l: norm of A}
  Let $A=(a_{i,j})$ be an $n \times n$ matrix whose entries are
  i.i.d subgaussian random variables. Then
  \[
     \P (\norm{A: B_2^n \to B_2^n} \ge \Cref{l: norm of A} \sqrt{n} )
     \le \exp (- n).
  \]
\end{lemma}

We will also need the volumetric estimate of the covering numbers
$N(K,D,t)$ (see e.g. \cite{P}). Denote by $|K|$ the volume of $K
\subset \R^n$.
\begin{lemma}  \label{l: volumetric}
Let $t>0$ and let $K,D \subset \R^n$ be convex symmetric bodies.
If $tD \subset K$, then
\[
  N(K,D,t)
  \le  \frac{3^n |K|}{|tD|}.
\]
\end{lemma}

\section{Hal\'{a}sz type lemma.}   \label{s: Halasz}

Let $\xi_1 \etc \xi_n$ be independent centered random variables.
To obtain the small ball probability estimates below, we have to
bound the probability that $\sum_{j=1}^n \xi_j$ is concentrated in
a small interval. One standard method of obtaining such bounds is
based on Berry-Ess\'{e}en Theorem (see, e.g. \cite{LPRT}).
However, this method has certain limitations. In particular, if
$\xi_j=t_j \e_j$, where $t_j \in [1,2]$ and $\e_j$ are $\pm1$
random variables, then Berry-Ess\'{e}en Theorem does not ``feel''
the distribution of the coefficients $t_j$, and thus does not
yield bounds better than $c/\sqrt{n}$ for the small ball
probability. To obtain better bounds we use the approach developed
by Hal\'{a}sz \cite{Ha1, Ha2}.

\begin{lemma}  \label{l: Halasz}
  Let $c>0, \ 0<\D <a/(2\pi)$ and let $\xi_1 \etc \xi_n$
  be independent
  random variables such that $\E \xi_i=0$, $\P (\xi_i>a) \ge c$ and
  $\P(\xi_i<-a) \ge c$.
  For $y \in \R$ set
  \[
    S_{\D}(y)=  \sumn \P(\xi_j -\xi_j' \in [y- \pi \D, y+ \pi \D]),
  \]
  where $\xi_j'$ is an independent copy of $\xi_j$.
  Then for any $v \in \R$
  \[
    \P \left ( \left |\sumn \xi_j -v \right | < \D \right )
    \le \frac{C}{n^{5/2} \D} \int_{3a/2}^{\infty} S_{\D}^2(y) \, dy
      +c e^{-c'n}.
  \]
\end{lemma}

\begin{proof}
 For
$t \in \R$ define
\begin{align*}
    \f_k(t) &= \E \exp (i \xi_k t)
    \intertext{and set}
    \f(t) &= \E \exp \left ( i t \sum_{k=1}^n \xi_k \right ) =
    \prod_{k=1}^n \f_k(t).
\end{align*}
Then by a Lemma of Ess\'{e}en \cite{Es}, for any $v \in \R$
\[
   Q= \P \left ( \left |\sumn \xi_j -v \right | < \D \right )
      \le c \int_{[-\pi/2,\pi/2]} |\f(t/\D)| \, dt.
\]
Let $\xi_k'$ be an independent copy of $\xi_k$ and let
$\nu_k=\xi_k-\xi_k'$. Then $\P(|\nu_k|> 2a) \ge 2c^2=\bar{c}$. We
have
\begin{equation}  \label{def phi_k}
  |\f_k(t)|^2 = \E \cos \nu_k t
\end{equation}
and
\[
  |\f(t)| \le \left ( \prod_{k=1}^n \exp \left (-1+|\f_k(t)|^2
    \right ) \right )^{1/2}
  = \exp \left ( -\frac{1}{2}\sum_{k=1}^n (1-|\f_k(t)|^2) \right ).
\]
Define a new random variable $\t_k$ by conditioning on
$|\nu_k|>2a$. For a Borel set $A \subset \R$ put
\[
  \P \left ( \t_k \in A \right ) =
    \frac{\P \left ( \nu_k \in A \setminus [-2a,2a]\right )}
         {\P \left ( |\nu_k|>2a \right )}.
\]

 Then by \eqref{def phi_k},
\[
  1- |\f_k(t)|^2 \ge \E (1- \cos \t_k t)
  \cdot \P \left ( |\nu_k|>2a \right )
  \ge \bar{c} \cdot \E (1- \cos \t_k t),
\]
so
\[
  |\f(t)| \le \exp (-c' f(t)),
\]
where
\[
  f(t) =\E \sum_{k=1}^n (1- \cos \t_k t).
\]
Let $T(m,r)= \{t \mid f(t/\D) \le m, \ |t| \le r \}$ and let
\[
  M = \max_{|t| \le \pi/2} f(t/\D).
\]

To estimate $M$ from below, notice that
\begin{align*}
  M&= \max_{|t| \le \pi/2} f(t/\D) \ge \frac{1}{\pi}
     \int_{-\pi/2}^{\pi/2} \E \sum_{k=1}^n (1- \cos (\t_k/\D) t) \, dt
     \\
   &= \E \sum_{k=1}^n
       \left (1- \frac{2}{\pi} \cdot \frac{\sin (\t_k/\D) \pi/2}{\t_k/\D}
       \right )
    \ge cn,
\end{align*}
since $|\t_k|/\D>2a/\D>4\pi$.

 To estimate the measure of $T(m, \pi/2)$ we use the argument of
\cite{Ha1}. For reader's convenience we present a complete proof.

\begin{lemma}  \label{l: measure}
  Let $0<m<M/4$. Then
  \[
    |T(m,\pi/2)| \le c \sqrt{\frac{m}{M}} \cdot |T(M/4, \pi)|.
  \]
\end{lemma}
\begin{proof}
Let $l=\sqrt{M/4m}$. Taking the integer part if necessary, we may
assume that $l$ is an integer. For $k \in \N$ set
\[
  S_k=\{\sum_{j=1}^k t_j \mid t_j \in T(m,\pi/2) \}.
\]
Note that $S_1=T(m, \pi/2)$.  Since
\[
  1-\cos\a=2 \sin^2(\a/2)
\]
and
\[
  \sin^2 \left (\sum_{j=1}^k \a_j \right )
  \le \left (\sum_{j=1}^k |\sin \a_j| \right )^2
  \le k \sum_{j=1}^k \sin^2 \a_j,
\]
we conclude that $S_k \subset T(k^2 m,k \pi/2)$. For $k \le l$ we
have $k^2 m <M$, so $(-\pi/2,\pi/2)\setminus T(k^2 m, k \pi/2)
\neq \emptyset$. For a Borel set $A$ denote $\mu(A)=|A \cap
[-\pi,\pi]|$. Now we shall prove by induction that for all $k \le
l$
\[
  \mu(S_k) \ge (k/2) \cdot \mu(S_1).
\]
Obviously, $\mu(S_2)=|S_2| \ge 2 \cdot |S_1|$, so this inequality
holds for $k=2$. Assume that $\mu(S_{k-1}) \ge (k-1)/2 \cdot
\mu(S_1)$. Note that the sets $S_k$ are closed. Let $v \in
(-\pi/2,\pi/2)$ be a boundary point of $S_k$. Such point exists
since $(-\pi/2,\pi/2)\setminus S_k \neq \emptyset$. Let
$\{v_j\}_{j=1}^{\infty}$ be a sequence of points in
$(-\pi/2,\pi/2)\setminus S_k$ converging to $v$. Then $(v_j - S_1)
\cap S_{k-1} = \emptyset$, so by continuity we have
\[
  \mu((v - S_1) \cap S_{k-1})=0.
\]
Since the set $S_1$ is symmetric, this implies
\[
  \mu((v + S_1) \cup S_{k-1})=\mu(v+S_1)+\mu(S_{k-1}).
\]
Both sets in the right hand side are contained in $S_{k+1}$ (to
see it for $S_{k-1}$ note that $0 \in S_2$). Since $v+S_1 \subset
[-\pi,\pi]$, the induction hypothesis implies
\[
  \mu(S_{k+1}) \ge \mu(v+S_1)+\mu(S_{k-1})
  \ge \mu(S_1)+ \frac{k-1}{2} \cdot \mu(S_1)
  =\frac{k+1}{2} \cdot \mu(S_1).
\]
Finally, since $S_l \cap [-\pi,\pi] \subset T(l^2 m, \pi)$, we get
\[
  |T(l^2 m, \pi)| \ge \frac{l}{2} \cdot |T(m,\pi/2)|.
\]
\end{proof}

We continue to prove Lemma \ref{l: Halasz}. Since
\begin{align*}
  Q &\le C \int_{[-\pi/2,\pi/2]} |\f(t/\D)| \, dt
     \le C \int_{[-\pi/2,\pi/2]} \exp(-c' f(t/\D)) \, dt \\
    &\le \bar{C}  \int_0^n |T(m, \pi/2)| e^{-c'm} \, dm,
\end{align*}
Lemma \ref{l: measure} implies
\begin{equation}  \label{hal: Q}
  Q  \le \frac{C'}{\sqrt{M}} \cdot |T(\frac{M}{4},\pi)|
    +ce^{-C'  M/16}.
    \le \frac{C'}{\sqrt{M}} \cdot |T(\frac{M}{4},\pi)|
    +ce^{-c'  n}.
\end{equation}
Here for $m>M/4$ we used a trivial estimate $|T(m,\pi/2)| \le
\pi$.

To complete the proof we have to estimate the measure of $T=T(M/4,
\pi)$ from above. For any $t \in T$ we have
\[
  g(t) = \sum_{k=1}^n \E \cos (\t_k t/\D)
  \ge n- M/4 \ge n/2.
\]
Let $w(x)= (1-|x|/\pi)\cdot \chi_{[-\pi,\pi]}(x)$ and put
$W=\hat{w}$. Then $W \ge 0$ and $W(t) \ge c$ for $|t| \le \pi$.
Hence by Parceval's equality,
\begin{align*}
    |T| &\le
      \left ( \frac{n}{2} \right )^{-2} \int_T |g(t)|^2 \,
         \le
      C \left  ( \frac{n}{2} \right )^{-2}
      \int_{\R} W^2(t) |g(t)|^2 \, dt \\
        &=
      \frac{C}{n^2} \int_{\R} \left |\E \sum_{k=1}^n w(\t_k/\D-y)
                        \right |^2 \, dy.
\end{align*}
Since $w \le \chi_{[-\pi,\pi]}$, the last expression does not
exceed
\begin{align*}
  &\frac{C}{n^2} \int_{\R}
      \left ( \sum_{k=1}^n \P (\frac{\t_k}{\D} \in [y-\pi, y+\pi])
      \right )^2 \, dy \\
  \le &\frac{C}{n^2 \D} \int_{\R}
             \left ( \sum_{k=1}^n \P (\t_k \in [z-\pi \D, z+\pi \D])
             \right )^2 \, dz.
\end{align*}
Since $\t_k \notin [-2a,2a]$ and $\pi \D < a/2$, we can integrate
only over $\R \setminus [-3a/2, 3a/2]$.

 Substituting this estimate into
\eqref{hal: Q}, we get
\[
  Q \le \frac{C}{n^{5/2}\D} \int_{\R \setminus [-3a/2,3a/2]}
             \left ( \sum_{k=1}^n \P (\t_k \in [z-\pi \D, z+\pi \D])
             \right )^2 \, dz +c e^{-c'n}.
\]
To finish the proof, recall that the variables $\t_k$ are
symmetric. This allows to change the integration set in the
previous inequality to $(3a/2, \infty)$. Moreover, if $z \in
(3a/2, \infty)$, then
\[
  \P (\t_k \in [z-\pi \D, z+\pi \D])
  \le \frac{1}{\bar{c}} \cdot \P (\nu_k \in [z-\pi \D, z+\pi \D]),
\]
so the random variables $\t_k$ can be replaced by
$\nu_k=\xi_k-\xi_k'$.

\end{proof}

\begin{remark}
A more delicate analysis shows that the term $c e^{-c'n}$ in the
formulation of Lemma \ref{l: Halasz} can always be eliminated.
However, we shall not prove it since this term does not affect the
results below.
\end{remark}

We shall apply Lemma \ref{l: Halasz} to weighted copies of the
same random variable. To formulate the result we have to introduce
a new notion.
\begin{definition}  \label{d: profile}
Let $x \in \R^m$. For $\D>0$ define the $\D$-profile of the vector
$x$ as a sequence $\{P_k(x,\D)\}_{k=1}^{\infty}$ such that
  \[
    P_k(x,\D)=|\{j \mid |x_j| \in (k \D, (k+1)\D] \}.
  \]
\end{definition}

\begin{theorem}   \label{t: weighted}
  Let $\b$ be a random variable such that $\E \b=0$ and
  $\P (\b>c) \ge c',\ \P (\b<-c) \ge c'$ for some $c,c'>0$.
  Let $\b_1 \ldots \b_m$ be independent copies of $\b$.
  Let $\D>0$ and let $(x_1 \ldots x_m) \in \R^m$ be a vector
  such $a<|x_j| <\overline{\Cref{t: weighted}} \, a$
  for some $a>0$. Then for any $\D<a/(2 \pi)$ and
  for any $v \in \R$
  \[
    \P \left ( \left | \sum_{j=1}^m \b_j x_j -v \right | <\D \right )
    \le \frac{\Cref{t: weighted}}{m^{5/2}} \sum_{k=1}^{\infty} P_k^2(x,\D).
  \]
\end{theorem}

\begin{proof}
We shall apply Lemma \ref{l: Halasz} to the random variables
$\xi_j=x_j \b_j$.

Let $\mathcal{M}(\R)$ be the set of all probability measures on
$\R$. Consider the function $F: \mathcal{M}(\R) \to \R_+$ defined
by
\[
  F(\mu)=\int_{3a/2}^{\infty} \tilde{S}_{\D}^2(y) \, dy,
\]
where
\[
  \tilde{S}_{\D}(y)
  =  \sum_{j=1}^m \mu(\frac{1}{|x_j|} \cdot [y- \pi \D, y+ \pi \D]).
\]

Since $F$ is a convex function on $\mathcal{M}(\R)$, it attains
the maximal value at an extreme point of this set, i.e. at some
delta-measure $\d_t, \ t \in \R$. Note that in this case
\[
  \tilde{S}_{\D}(y)= |\{j \mid t|x_j| \in [y- \pi \D, y+ \pi \D]
                      \}
  =\sum_{j=1}^m \chi(t|x_j|-y),
\]
where $\chi=\chi_{[-\pi \D,\pi \D]}$ is the indicator function of
$[-\pi \D,\pi \D]$. For $t<\frac{1}{2C}$ we have $t|x_j|<a/2$, so
$\tilde{S}_{\D}(y)=0$ for any $y \ge 3a/2$, and thus $F(\d_t)=0$.
If $t \ge \frac{1}{2C}$, then
\begin{align*}
    F(\d_t)
  &=  \sum_{j=1}^m \sum_{l=1}^m
  \int_{3a/2}^{\infty} \chi(t|x_j|-y) \chi(t|x_l|-y) \, dy \\
  &\le 2 \pi \D |\{ (j,l) \mid t\big | |x_j|-|x_l| \big | \le \pi
  \D \}|=g(t).
\end{align*}
Since the function $g$ is decreasing,
\begin{align*}
  F(\d_t)
  &\le g(\frac{1}{2C}) \le 2 \pi \D
  \sum_{l=1}^{\infty}
  |\{j \mid \Big | |x_j|- l\D \Big | \le 2\pi \D \cdot C\}|^2 \\
  &\le \bar{C} \D \sum_{k=1}^{\infty}
  |\{j \mid |x_j| \in (k \D, (k+1)\D] \}|^2.
\end{align*}
The last inequality holds since we can cover each interval $[l
\D-2 \pi \D C, l \D +2 \pi \D C]$ by at most $2 \pi C +2$
intervals $(k \D, (k+1)\D]$.

Let $\mu$ be the distribution of the random variable $\b-\b'$,
where $\b'$ is an independent copy of $\b$. Applying Lemma \ref{l:
Halasz} to the random variables $\xi_j=x_j \cdot \b_j$, we have
\begin{align*}
  \P &\left ( \left |\sum_{j=1}^m \b_j x_j -v \right | < \D \right )
    \le \frac{C}{m^{5/2} \D} F(\mu) +c e^{-c'm}\\
    &\le \frac{C'}{m^{5/2}} \sum_{k=1}^{\infty}
  |\{j \mid |x_j| \in (k \D, (k+1)\D] \}|^2 +c e^{-c'm}.
\end{align*}
Since the sum in the right hand side is at least $m$, the second
term is negligible compare to the first one. Thus,
\[
  \P \left ( \left |\sum_{j=1}^m \b_j x_j -v \right | < \D \right )
    \le \frac{2C'}{m^{5/2}} \sum_{k=1}^{\infty}
  |\{j \mid |x_j| \in (k \D, (k+1)\D] \}|^2.
\]
\end{proof}

\section{Small ball probability estimates.}   \label{s: small ball}

Let $G$ be an $n \times n$ Gaussian matrix. If $x \in S^{n-1}$ is
any unit vector, then $y=Gx$ is the standard Gaussian vector in
$\R^n$. Hence for any $t>0$ we have $\P(|y_j|<t) \le t \cdot
\sqrt{2/\pi}$ for any coordinate. Moreover,
\[
  \P(\norm{y} \le t \cdot \sqrt{n})
  \le (2 \pi)^{-n/2} \vol (t \sqrt{n} B_2^n)
  \le (Ct)^n.
\]
We would like to have the same small ball probability estimates
for the random vector $y=Ax$. However, it is easy to see that it
is impossible to achieve such estimate for {\em all} directions $x
\in S^{n-1}$. Indeed, if $A$ is a random $\pm 1$ matrix and
$x=(1/\sqrt{2}, 1/\sqrt{2},0 \ldots 0)$, then $\P(y_j=0)=1/2$ and
$\P(y=0)=2^{-n}$. Analyzing this example, we see that the reason
that the small ball estimate fails is the concentration of the
Euclidean norm of $x$ on a few coordinates. If the vector $x$ is
``spread'', we can expect a more regular behavior of the small
ball probability.

Although we cannot prove the Gaussian type estimates for all
directions and all $t>0$, it is possible to obtain such estimates
for {\em most} directions provided that $t$ is sufficiently large
($t>t_0)$. Moreover, the more we assume about the regularity of
distribution of the coordinates of $x$, the smaller value of $t_0$
we can take. This general statement is illustrated by the series
of results below.

The first result is valid for any direction. The following Lemma
 is a particular case of
\cite{LPRT}, Proposition 3.4.
\begin{lemma}  \label{l: peaked}
  Let $A$ be an $n \times n$ matrix with i.i.d. subgaussian
  entries. Then for every $x \in S^{n-1}$
  \[
    \P(\norm{Ax} \le \Cref{l: peaked} \sqrt{n})
    \le \exp (-\cref{l: peaked}n).
  \]
\end{lemma}
The example considered at the beginning of this section shows that
this estimate cannot be improved for a general random matrix.

If we assume that all coordinates of the vector $x$ are
comparable, then we have the following Lemma, which is a
particular case of Proposition 3.4 \cite{LPRTV2} (see also
Proposition 3.2 \cite{LPRT}).
\begin{lemma}  \label{l: BE}
  Let $\b$ be a random variable such that $\E \b=0, \ \E \b^2=1$ and
  let $\b_1 \etc \b_m$ be independent copies of $\b$. Let $0<r<R$ and let $x_1
  \etc x_m \in \R$ be such that $r/\sqrt{m} \le |x_j| \le
  R/\sqrt{m}$ for any $j$. Then for any
  $t \ge \cref{l: BE}/\sqrt{m}$
  and for any $v \in \R$
  \[
    \P \left ( \left | \sum_{j=1}^m \b_j x_j -v \right | <t
       \right )
    \le \Cref{l: BE} t.
  \]
  Here $\cref{l: BE}$ and $\Cref{l: BE}$ depend only on $r$ and
  $R$.
\end{lemma}
\begin{proof}
The proof is based  on  Berry-Ess\'een theorem (cf., e.g.,
\cite{Strook}, Section 2.1).

\begin{theorem}
\label{t: be} Let $(\zeta_j)_{i=1}^m$ be a sequence of independent
random variables with expectation $0$ and  finite third moments,
and let $A^2 : = \sum^m_{j=1}\E |\zeta_j|^2 $.  Then for every
$\tau \in \R$ one has
$$
\Bigl| \P\Bigl( \sum^m_{j=1} \zeta_j < \tau A \Bigr) - \P \left( g
< \tau \right ) \Bigr| \le  (c/A^3) \sum^m_{j=1} \E |\zeta_j|^3,
$$
where $g$ is a Gaussian random variable with $N(0,1)$ distribution
and $c \geq 1$ is a universal constant.
\end{theorem}

Let $\zeta_j = \b_j x_j$. Then $A^2 := \sum_{j=1}^m  \E \zeta_j^2
= \norm{x}^2 \ge r^2$. Since the random variables $\b_j$ are
copies of a subgaussian random variable $\b$, $\E |\b|^3 \le C$
for some absolute constant $C$. Hence, $\E \sum_{j=1}^m
|\zeta_j|^3 \le C \sum_{j=1}^m |x_j|^3 \le C'/\sqrt{m}$. By
Theorem~\ref{t: be} we get
\begin{align*}
  \P \left ( \left | \sum_{j=1}^m \b_j x_j -v \right | <t \right )
  &\le \P \left ( \frac{v-t}{c} \le g  < \frac{v+t}{c} \right )
   + \frac{c'}{\sqrt{m}}  \\
  &\le C''t + \frac{c'}{\sqrt{m}}
   \le 2C''t,
\end{align*}
provided $t \ge \frac{C''}{c' \sqrt{m}}$.
\end{proof}

If $x=(1/\sqrt{m}, \ldots, 1/\sqrt{m})$, then
  \[
    \P \left ( \left | \sum_{j=1}^m \b_j x_j \right | =0
       \right )
    \ge C/\sqrt{m}.
  \]
This shows that the bound $t \ge \cref{l: BE}/\sqrt{m}$ in Lemma
\ref{l: BE} is necessary.

The proofs of Lemma \ref{l: peaked} and Lemma \ref{l: BE} are
based on Paley--Zygmund inequality and Berry--Ess\'{e}en Theorem
respectively. To obtain the linear decay of small ball probability
for $t \le \cref{l: BE}/\sqrt{m}$, we use the third technique,
namely Hal\'{a}sz method. However, since the formulation of the
result requires several technical assumptions on the vector $x$,
we postpone it to Section \ref{s: regular profile}, where these
assumptions appear.

To translate the small ball probability estimate for a single
coordinate to a similar estimate for the norm we use the Laplace
transform technique, developed in \cite{LPRT}. The following Lemma
improves the argument used in the proof of Theorem 3.1
\cite{LPRT}.
\begin{lemma} \label{l: Laplace transform}
  Let $\D>0$ and let $Y$ be a random variable such that for any $v
  \in \R$ and for any $t\ge \D, \ \P (|Y-v|>t) \le Lt$. Let
  $y=(Y_1 \etc Y_n)$ be a random vector, whose coordinates are
  independent copies of $Y$. Then for any $z \in \R^n$
  \[
    \P \left ( \norm{y-z} \le \D \sqrt{n} \right )
    \le (\Cref{l: Laplace transform} L \D)^n.
  \]
\end{lemma}
\begin{proof}
We have
\begin{align*}
   \P \left ( \norm{y-z} \le  \D \sqrt{n} \right )
   &= \P \left ( \sum_{i=1}^n (Y_i-z_i)^2 \le  \D n
         \right ) \\
   &= \P \left (  n - \frac{1}{ \D} \sum_{i=1}^n (Y_i-z_i)^2 \ge 0
         \right ) \\
   &\le \E \exp \left ( n - \frac{1}{ \D} \sum_{i=1}^n (Y_i-z_i)^2
                \right ) \\
   &= e^n \cdot \prod_{i=1}^n \E \exp (- \frac{1}{ \D} (Y_i-z_i)^2).
\end{align*}
To estimate the last expectation we use Lemma \ref{l: regular
disrepancy}.
\begin{align*}
  \E \exp (- \frac{1}{ \D} (Y_i-z_i)^2)
  &= \int_0^1 \P \left (\exp \Big ( - \frac{1}{ \D} (Y_i-z_i)^2 \Big )
                     >s
                \right ) \, ds \\
  &= \int_0^{\infty} 2u e^{-u^2} \P (|Y_i-z_i| <  \D u) \, du \\
  &\le \int_0^{1} u e^{-u^2/2} L \D \, du \\
  &+     \int_{1}^{\infty} u e^{-u^2/2} L \D u \, du \\
  &\le \bar{C} L \D.
\end{align*}
Substituting this into the previous inequality, we get
\[
   \P \left ( \norm{y-z} \le  \D \sqrt{n} \right )
   \le (e \cdot \bar{C} L \D)^n.
\]
\end{proof}

\section{Partition of the sphere.}   \label{s: sphere}

To apply the small ball probability estimates proved in the
previous section we have to decompose the sphere into different
regions depending on the distribution of the coordinates of a
point. We start by decomposing the sphere $S^{n-1}$ in two parts
following \cite{LPRT, LPRTV1, LPRTV2}. We shall define two sets:
$V_P$ -- the set of vectors, whose Euclidean norm is concentrated
on a few coordinates, and $V_S$ -- the set of vectors whose
coordinates are evenly spread. Let  $r<1<R$ be the numbers to be
chosen later. Given $x=(x_1 \etc x_n) \in S^{n-1}$, set $\s(x)=\{i
\mid |x_i| \le R/ \sqrt{n} \}$. Let $P_I$ be the coordinate
projection on the set $I \subset \{ \etc n\}$. Set
\begin{align*}
    V_P &= \{x \in S^{n-1} \mid \norm{P_{\s(x)}x} <  r \} \\
    V_S &= \{x \in S^{n-1} \mid \norm{P_{\s(x)}x} \ge r \}.
\end{align*}
First we shall show that with high probability $\norm{Ax} \ge C
\sqrt{n}$ for any $x \in V_P$.

For a single vector $x \in \R^n$ this probability was estimated in
Lemma \ref{l: peaked}. We shall combine this estimate with an
$\e$-net argument.
\begin{lemma}   \label{l: entropy of V_P}
  For any $r<1/2$
  \[
    \log N(V_P, B_2^n,2r) \le
    \frac{n}{R} \cdot \log \left (\frac{3R}{r} \right ).
  \]
\end{lemma}

\begin{proof}
If $x \in B_2^n$, then $|\{1 \etc n\} \setminus \s(x)| \le n/R$.
Hence, the set $V_P$ is contained in the sum of two sets: $rB_2^n$
and
\[
  W_P=\{x \in B_2^n \mid |\supp (x)| \le n/R^2\}.
\]
Since $W_P$ is contained in the union of unit balls in all
coordinate subspaces of dimension $l=n/R$, Lemma \ref{l:
volumetric} implies
\[
  N(W_P, B_2^n, r)
  \le \binom{n}{l} \cdot
  N(B_2^l, B_2^{l}, r)
  \le \binom{n}{l} \cdot \left (\frac{3}{r} \right )^l.
\]
Finally,
\[
  \log N(V_P,B_2^n,2r) \le \log N(W_P,B_2^n,r)
  \le l \cdot \log \left (\frac{3n}{lr} \right )
  \le \frac{n}{R} \cdot \log \left (\frac{3R}{r} \right ).
\]
\end{proof}

Recall that $\Cref{l: peaked} < \Cref{l: norm of A}$. Set
$r=\Cref{l: peaked}/2\Cref{l: norm of A}$ and choose the number
$R>1$ so that
\[
  \frac{1}{R} \cdot \log \left (\frac{3R}{r} \right )
    < \frac{\cref{l: peaked}}{2}.
\]
For these parameters we prove that the norm of $Ax$ is bounded
below for all $x \in V_P$ with high probability.
\begin{lemma}  \label{l: all peaked}
  \[
    \P \left ( \exists x \in V_P \mid \norm{Ax}
               \le \Cref{l: peaked} \sqrt{n}/2
       \right )
    \le 2 \exp (-\cref{l: peaked} n).
  \]
\end{lemma}

\begin{proof}
By Lemma \ref{l: entropy of V_P}, the set $V_P$ contains a
$(\Cref{l: peaked}/2\Cref{l: norm of A})$-net $\mathcal{N}$ in the
$\ell_2$-metric of cardinality at most $\exp (\cref{l: peaked}
n/2)$. Let
\begin{align*}
  \Om_0&= \{\om \mid \norm{A}>\Cref{l: norm of A} \sqrt{n} \}
  \intertext{and let}
  \Om_P&= \{\om \mid
                \exists x \in \mathcal{N} \
                \norm{A(\om)x} \le \Cref{l: peaked} \sqrt{n}
         \}.
\end{align*}
Then Lemma \ref{l: peaked} implies
\[
  \P ( \Om_0 )+\P(\Om_P)
  \le \exp(-n) + \exp (- \cref{l: peaked}n)
  \le 2\exp(- \cref{l: peaked}n).
\]
Let $\om \notin \Om_P$. Pick any $x \in V_P$. There exists $y \in
\mathcal{N}$ such that $\norm{x-y}_2 \le \Cref{l:
peaked}/2\Cref{l: norm of A}$. Hence
\begin{align*}
  \norm{Ax}
  &\ge \norm{Ay} - \norm{A(x-y)}
  \ge \Cref{l: peaked} \sqrt{n} - \norm{A: B_2^n \to B_2^n}  \cdot
       \norm{x-y}_2 \\
  &\ge \frac{\Cref{l: peaked}}{2} \sqrt{n}.
\end{align*}

\end{proof}

For $x=(x_1 \etc x_n) \in V_S$ denote
\begin{equation}  \label{def J}
  J(x) = \left \{j \mid \frac{r}{2\sqrt{n}} \le |x_j|
                        \le \frac{R}{\sqrt{n}} \right \}.
\end{equation}
Note that
\[
 \sum_{j \in J(X)} x_j^2
 \ge \sum_{j \in \s(X)} x_j^2 - \frac{r^2}{2}
 \ge \frac{r^2}{2},
\]
 so
\[
  |J(x)| \ge (r^2/2 R^2) \cdot n=:m.
\]

Let $0<\D< r/2\sqrt{n}$ be a number to be chosen later. We shall
cover the interval $[\frac{r}{2 \sqrt{n}}, \frac{R}{\sqrt{n}}]$ by
\[
  k=\left \lceil\frac{R-r/2}{\sqrt{n}\D} \right \rceil
\]
    consecutive intervals $(j\D,(j+1)\D]$,
where $j=k_0, (k_0+1) \etc (k_0+k)$, and $k_0$ is the largest
number such that $k_0 \D < r/2\sqrt{n}$.
 Then we shall decompose the set $V_S$
in two subsets: one containing the points whose coordinates are
concentrated in a few such intervals, and the other containing
points with evenly spread coordinates. This will be done using the
$\D$-profile, defined in \ref{d: profile}. Note that if $m$
coordinates of the vector $x$ are evenly spread among $k$
intervals, then
\[
  \sum_{i=1}^{\infty} P_i^2(x,\D) \sim \frac{m^2}{k}
  \sim m^{5/2} \D.
\]
This observation leads to the following

\begin{definition} \label{d: regular profile}
  Let $\D>0$ and let $Q>1$.
  We say that a vector $x\in V_S$ has a $(\D,Q)$-regular profile
  if there
  exists a set $J \subset J(x)$ such that $|J| \ge m/2$ and
  \[
    \sum_{i=1}^{\infty} P_i^2(x|_J,\D)
    \le Q m^{5/2} \D
    =: \Cref{d: regular profile} Q \cdot \frac{m^2}{k}.
  \]
  Here $x |_J \in \R^n$ is a vector with coordinates
  $x|_J(j)=x(j) \cdot \chi_J(j)$.

  If such set $J$ does not exist, we call $x$ a vector of
  $(\D,Q)$-singular
  profile.
\end{definition}

Note that $\sum_{i=1}^{\infty} P_i^2(x|_J,\D) \ge m/2$. Hence, if
$\D<m^{-3/2}/2$, then every vector in $V_S$ will be a vector of a
$(\D,Q)$-singular profile.

Vectors of regular and singular profile will be treated
differently. Namely, in Section \ref{s: regular profile} we prove
that vectors of regular profile satisfy the small ball probability
estimate of the type $Ct$ for $t \ge \D$. This allows to use
conditioning to estimate the probability that $\norm{Ax}$ is small
for some vector $x$ of regular profile. In Section \ref{s:
singular profile} we prove that the set of vectors of singular
profile admits a small $\e$-net. This fact combined with Lemma
\ref{l: BE} allows to estimate the probability that there exists a
vector $x$ of singular profile such that $\norm{Ax}$ is small
using the standard $\e$-net argument.

\section{Vectors of a regular profile.}   \label{s: regular profile}

To estimate the small ball probability for a vector of a regular
profile we apply Theorem \ref{t: weighted}.

\begin{lemma}  \label{l: regular disrepancy}
  Let $\D \le \frac{r}{4 \pi \sqrt{n}}$.
  Let $x \in V_S$ be a vector of $(\D,Q)$-regular profile.
  Then for any $t \ge \D$
  \[
     \P \left ( \left | \sum_{j=1}^n \b_j x_j -v \right | <t \right )
     \le \Cref{l: regular disrepancy} Q \cdot t.
  \]
\end{lemma}
\begin{proof}
Let $J \subset \{1 \etc n\}, \ |J| \ge m/2$ be the set from
Definition \ref{d: regular profile}. Denote by $\E_{J^c}$ the
expectation with respect to the random variables $\b_j$, where $j
\in J^c=\{1 \etc n\} \setminus J$. Then
\begin{align*}
  &\P \left ( \left | \sum_{j=1}^n \b_j x_j -v \right | <t
      \right ) \\
     &= \E_{J^c}
     \P \left ( \left | \sum_{j \in J} \b_j x_j -(v
               + \sum_{j \in J^c} \b_j x_j)
               \right | <t
               \mid \b_j, \ j \in J^c
       \right )
\end{align*}
Hence, it is enough to estimate the conditional probability.

Recall that $\b$ is a centered subgaussian random variable of
variance $1$. It is well-known that such variable satisfies
$\P(\b>c) \ge c', \ \P(\b<-c) \ge c'$ for some absolute constants
$c,c'$. Moreover, a simple Paley--Zygmund type argument shows that
this estimates hold if we assume only that $\E \b =0$ and the
second and the fourth moment of $\b$ are comparable. Hence, for
$t=\D$ the Lemma follows from Theorem \ref{t: weighted}, where we
set $a=r/\sqrt{n}, \ \overline{\Cref{t: weighted}}=R/r$.

To prove the Lemma for other values of $t$, assume first that
$t=\D_s =2^s \D< \frac{r}{4 \pi \sqrt{n}}$ for some $s \in \N$.
Consider the $\D_s$-profile  of $x|_J$:
\[
  P_l(x|_J, \D_s)= |\{j \in J \mid |x_j| \in (l \D_s, (l+1)\D_s] \}|.
\]
Notice that each interval $(l \D_s, (l+1)\D_s]$ is a union of
$2^s$ intervals $(i \D, (i+1)\D]$. Hence
\[
  \sum_{l=1}^{\infty} {P}_l^2(x|_J,\D_s)
  \le 2^s
  \sum_{i=1}^{\infty} P_i^2(x|_J,\D)
  \le 2^s Q m^{5/2} \D
  = Q m^{5/2} t.
\]
Applying Theorem \ref{t: weighted} with $\D$ replaced by $\D_s$
and $v'=v + \sum_{j \in J^c} \b_j x_j$, we obtain
\[
   \P \left ( \left | \sum_{j \in J} \b_j x_j -(v
               + \sum_{j \in J^c} \b_j x_j)
               \right | <t
               \mid \b_j, \ j \in J^c
       \right ) \le \Cref{t: weighted} Q t.
\]
For $2^s\D <t <2^{s+1}\D$ the result follows from the previous
inequality applied for $t=2^s \D$. If $t \ge \cref{l:
BE}/\sqrt{m}=\frac{\sqrt{2}\, \cref{l: BE}  R}{r \sqrt{n}}$, Lema
\ref{l: BE} implies
\[
   \P \left ( \left | \sum_{j \in J} \b_j x_j -(v
               + \sum_{j \in J^c} \b_j x_j)
               \right | <t
               \mid \b_j, \ j \in J^c
       \right ) \le \Cref{l: BE} t \le \Cref{l: BE} Q t.
\]
Finally, if $\frac{r}{4 \pi \sqrt{n}} < t < \frac{\sqrt{2}\,
\cref{l: BE}  R}{r \sqrt{n}}$, the previous inequality applied to
$t_0=\frac{\sqrt{2}\, \cref{l: BE}  R}{r \sqrt{n}}$ implies
\[
   \P \left ( \left | \sum_{j \in J} \b_j x_j -(v
               + \sum_{j \in J^c} \b_j x_j)
               \right | <t
               \mid \b_j, \ j \in J^c
       \right ) \le \Cref{l: BE} Q t_0 \le CQ t,
\]
where $C=\Cref{l: BE} \cdot \frac{\sqrt{2}\, \cref{l: BE}  R}{r }
\cdot \frac{4 \pi}{r}$.
\end{proof}

Now we  estimate the probability that $\norm{A(\om)x}$ is small
for some vector of a regular profile.

\begin{theorem}  \label{t: regular profile}
  Let $\D>0$ and let $U$ be the set of vectors of
  $(\D,Q)$-regular profile. Then
  \[
    \P \left (\exists \, x \in U
              \mid \norm{Ax} \le \frac{\D}{2 \sqrt{n}}
       \right )
    \le \Cref{l: regular disrepancy} Q \D n.
  \]
\end{theorem}

\begin{proof}
Set
\[
  s=\frac{\D}{2 \sqrt{n}}.
\]
Let $\Om$ be the event described in Theorem \ref{t: regular
profile}.
  Denote the rows of $A$ by $a_1 \etc a_n$.
 Note that since
$\norm{A^{-1}}=\norm{(A^{-1})^T}$, for any $\om \in \Om$ there
exists a vector $u=(u_1 \etc u_n) \in S^{n-1}$ such that
\[
  u_1 a_1  + \ldots + u_n a_n =z,
\]
where $\norm{z}<s$. Then $\Om = \cup_{k=1}^n \Om_k$, where $\Om_k$
is the event $|u_k| \ge 1/\sqrt{n}$. Since the events $\Om_k$ have
the same probability, it is enough to estimate $\P (\Om_n)$.

To this end we condition on the first $n-1$ rows of the matrix
$A=A(\om)$:
\[
  \P(\Om_n)=\E_{a_1 \etc a_{n-1}} \P (\Om_n \mid a_1 \etc
  a_{n-1}).
\]
Here $\E_{a_1 \etc a_{n-1}}$ is the expectation with respect to
the first $n-1$ rows of the matrix $A$. Take {\em any} vector $y
\in U$ such that
\[
  \sum_{j=1}^{n-1} \pr{a_j}{y}^2 < s^2.
\]
If such vector does not exist, then $\norm{Ay} \ge s$ for all $y
\in U$, and so $\om \notin \Om$. Note that the vector $y$ can be
chosen using only  $a_1 \etc a_{n-1}$. We have
\[
  a_n = \frac{1}{u_n} (u_1 a_1 + \ldots + u_{n-1} a_{n-1} -z),
\]
so for $\om \in \Om_n$
\begin{align*}
  |\pr{a_n}{y}|
  &= \frac{1}{|u_n|}
     \left | \sum_{j=1}^{n-1} u_j \pr{a_j}{y} - \pr{z}{y}
     \right | \\
  &\le \sqrt{n} \left ( \left ( \sum_{j=1}^{n-1} u_j^2
                        \right )^{1/2}
                        \left ( \sum_{j=1}^{n-1} \pr{a_j}{y}^2
                        \right )^{1/2} + \norm{z}
                \right )
   \le 2 \sqrt{n} \cdot s= \D.
\end{align*}

 The row $a_n$ is independent of $a_1 \etc a_{n-1}$.
 Hence, Lemma \ref{l: regular disrepancy} implies
\begin{align*}
  \P (\Om_n)
  &\le \P \left (|\pr{a_n}{y}| \le \D
                \mid a_1 \etc a_{n-1}\right ) \\
  &= \P \left ( \left | \sum_{j=1}^n \b_{n,j} y_j \right | \le \D
                 \mid a_1 \etc a_{n-1}\right )
  \le \Cref{l: regular disrepancy} Q \D,
\end{align*}
and so $\P(\Om) \le \Cref{l: regular disrepancy} Q \D n$.
\end{proof}

\section{Vectors of a singular profile.}   \label{s: singular profile}

We prove first that the set of vectors of singular profile admits
a small $\D$-net in the $\ell_{\infty}$-metric.

\begin{lemma}  \label{l: net singular}
  Let $ \overline{\Cref{l: net singular}} n^{-3/2}
  \le \D \le n^{-1/2}$, where
  $\overline{\Cref{l: net singular}}= \frac{2R^3}{r^2}$ and
   let $W_S$ be the set of vectors of $(\D,Q)$-singular profile.
  Let $\eta>0$ be such that
  \[
    C(\eta) < \Cref{d: regular profile} Q,
  \]
  where $C(\eta)$ is the function defined in Lemma \ref{l:
  concentration}.
   Then there
  exists a $\D$-net $\mathcal{N}$ in $W_S$ in $\ell_{\infty}$-metric
  such that
  \[
     |\mathcal{N}|
     \le \left ( \frac{\Cref{l: net singular}}{\D \sqrt{n}}
                 \eta^{\cref{l: net singular}}\right )^n.
  \]
\end{lemma}
\begin{remark}
Lemma \ref{l: volumetric} implies that there exists a $\D$-net for
$S^{n-1}$ in the  $\ell_{\infty}$-metric with less than $(C \D/
\sqrt{n})^n$ points. Thus, considering only vectors of a singular
profile, we gain the factor $\eta^{\cref{l: net singular} \cdot
n}$ in the estimate of the size of a $\D$-net.
\end{remark}

\begin{proof}
Let $J \subset \{1 \etc n\}$ and denote $J'=\{1 \etc n\} \setminus
J$. Let $W_J \subset W_S$ be the set of all vectors $x$ of a
$(\D,Q)$-singular profile for which $J(x)=J$. We shall construct
$\D$-nets in each $W_J$ separately. To this end we shall use Lemma
\ref{l: concentration} to construct a $\D$-net for the set $P_J
W_J$, where $P_J$ is the coordinate projection on $\R^J$. Then the
product of this $\D$-net and a $\D$-net for the ball $B_2^{J'}$
will form a $\D$-net for the whole $W_J$.

 Assume that $J=\{1
\etc l\}$, where $l \ge m$. Let $I_1 \etc I_k$ be consecutive
subintervals $(i \D, (i+1) \D],\ i=k_0 \etc k_0+k$, covering the
interval $[\frac{r}{2 \sqrt{n}}, \frac{R}{\sqrt{n}}]$, which
appear in the definition of profile. Recall that
\[
  k=\left \lceil\frac{R-r/2}{\sqrt{n}\D} \right \rceil
\]
 The restriction on $\D$ implies that $k \le m$. Let $d_i$
be the center of the interval $I_i$. Set
\[
  \mathcal{M}_J=\{x \in \R^J \mid
                |x_j| \in \{d_1 \etc d_k \} \
               \text{for } j\in J \}.
\]
Then $|\mathcal{M}_J| = (2k)^l$. Let $\mathcal{N}_J$ be the set of
all $x \in \mathcal{M}_J$ for which there exists a vector $y \in
W_J$ such that $-\D/2<y_j-x_j \le \D/2$ for all $j \in J$. The set
$\mathcal{N}_J$ forms a $\D$-net for $W_J$ in the $\ell_{\infty}$
metric. To estimate its cardinality we use the probabilistic
method.

Let  $X(1) \etc X(l)$ be independent random variables uniformly
distributed on the set $\{1 \etc k\}$. Let $N \subset \{1 \etc
k\}^l$ be the set of all $l$-tuples $(v(1) \etc v(l))$ such that
$|x_j|=d_{v(j)}, \ j=1 \etc l$ for some $x=(x_1 \etc x_l) \in
\mathcal{N}_J$. Since both $\mathcal{M}_J$ and $\mathcal{N}_J$ are
invariant under changes of signs of the coordinates,
\[
  \P((X(1) \etc X(l)) \in N)=
  \frac{|\mathcal{N}_J|}{|\mathcal{M}_J|}.
\]
Let $(X(1) \etc X(l)) \in N$ and let $x \in \R^l$ be such that
$x_j = d_{X(j)}$. Let $y \in W_J$ be a vector such that
$-\D/2<y_j-x_j \le \D/2$ for all $j \in J$. Then for any $j \in J,
\ y_j \in I_i$ implies that $X(j)=i$. Let $E \subset J$ be any set
containing at least $m/2$ elements. Then
\[
  \sum_{i=1}^{\infty} P_i^2(y|_E,\D)
  = \sum_{i=1}^k |\{j \in E \mid X(j)=i\}|^2.
\]
Since $y$ is a vector of a singular profile, this implies
\[
  \sum_{i=1}^k |\{j \in E \mid X(j)=i\}|^2
  \ge  Q m^{5/2} \D
  = \Cref{d: regular profile} \cdot Q \frac{m^2}{k}
  >C(\eta) \cdot \frac{m^2}{k}.
\]

 Now Lemma \ref{l:
concentration} implies that $\P(N) \le \eta^{l}$, so
\[
  |\mathcal{N}_J| \le (2k \eta)^l
  = \left ( \frac{R-2r}{\D \sqrt{n}}\eta \right )^l.
\]
To estimate the cardinality of the $\D$-net for the whole $W_J$ we
use Lemma \ref{l: volumetric}. Since $\D \le 1/\sqrt{|J|}, \ \D
B_{\infty}^J \subset B_2^J$, so
\[
   N(P_{J'} W_J, B_{\infty}^{J'}, \D)
   \le N(B_2^{J'}, B_{\infty}^{J'}, \D)
   \le 3^{n-l} \frac{|B_2^{J'}|}{|\D B_{\infty}^{J'}|}
   \le \left ( \frac{c}{\D \sqrt{n-l}} \right )^{n-l}.
\]
Since the function $f(t)=(a/t)^t$ is increasing for $0<t<a/e$, the
right-hand side of the previous inequality is bounded by $\left (
c / \D \sqrt{n} \right )^{n}$. Hence,
\begin{align*}
  N(W_J, B_{\infty}^n, \D)
  &\le N(P_J W_J, B_{\infty}^J, \D)
      \cdot N(P_{J'} W_J, B_{\infty}^{J'}, \D) \\
  &\le |\mathcal{N}_J| \cdot
      \left ( \frac{c}{\D \sqrt{n}} \right )^{n}
   \le \left ( \frac{c'}{\D \sqrt{n}} \eta^{l/n} \right )^n
\end{align*}

  Finally, set
\[
  \mathcal{N}= \bigcup_{|J| \ge m} \mathcal{N}_J.
\]
Then
\[
  |\mathcal{N}| \le \sum_{l=m}^n \sum_{|J|=l} |\mathcal{N}_J|
  \le 2^n \left ( \frac{c'}{\D \sqrt{n}} \eta^{m/n} \right )^n.
\]
Thus, Lemma \ref{l: net singular} holds with $\cref{l: net
singular}=m/n=\frac{r^2}{2R^2}$.
\end{proof}

Now we are ready to show that $\norm{Ax} \ge c$ for {\em all}
vectors of a $(\D,Q)$-singular profile with probability
exponentially close to 1.

\begin{theorem}  \label{t: singular profile}
  There exists an absolute constant $Q_0$ with the following
  property.
  Let $\D \ge \Cref{t: singular profile} n^{-3/2}$, where
  $\Cref{t: singular profile}=
  \max(\cref{l: BE},\overline{\Cref{l: net singular}})$.
   Denote by $\Om_{\D}$ the event
  that there exists a vector $x \in V_S$ of $(\D,Q_0)$-singular
  profile
  such that
  $\norm{Ax} \le \frac{\D}{2} n$.
  Then
  \[
    \P (\Om_{\D}) \le 3\exp(-n).
  \]
\end{theorem}

\begin{proof}
We consider two cases. First, we assume that $\D \ge \D_1=\cref{l:
BE}/ n$. In this case we estimate the small ball probability using
Lemma \ref{l: BE} and the size of the $\e$-net using Lemma \ref{l:
net singular}. Note that only the second estimate uses the  the
profile of the vectors. Then we conclude the proof with the
standard approximation argument.

The case $\D \le \D_1$ is more involved. From Case 1 we know that
there exists $Q_1$ such that  all vectors of $(\D_1,
Q_1)$-singular profile satisfy $\norm{Ax} \ge \frac{\D_1}{2} n$
with probability at least $1-e^{-n}$. Hence, it is enough to
consider only vectors whose profile is regular on the scale $\D_1$
and singular on the scale $\D$. For these vectors we use the
regular profile in  Lemma \ref{l: regular disrepancy} to estimate
the small ball probability and singular profile in Lemma \ref{l:
net singular} to estimate the size of the $\e$-net. The same
approximation argument finishes the proof.

{\bf Case 1}. Assume first that $\D \ge \D_1= \cref{l: BE}/ n$.
 Let $Q_1>1$
be a number to be chosen later. Let $\mathcal{M}$ be the smallest
$\frac{\D}{2 \Cref{l: norm of A}}$-net in the set of the vectors
of $(\D, Q_1)$-singular profile in $\ell_{\infty}$ metric.

Let $x \in V_S$ and let $J=J(x)$ defined in \eqref{def J}. Denote
$J^c= \{1 \etc n\} \setminus J$. Then Lemma \ref{l: BE} implies
\begin{align*}
  \P \left ( \left | \sum_{j=1}^n \b_j x_j \right | \le  t
     \right )
  &= \E_{J^c} \P \left ( \left | \sum_{j \in J} \b_j x_j
          + \sum_{j \in J^c} \b_j x_j \right | \le  t
          \mid \b_j, \ j \in J^c
     \right )  \\
  &\le \Cref{l: BE} t
\end{align*}
for all $t \ge \cref{l: BE} /\sqrt{n}$. Since $\D \sqrt{n} \ge
\cref{l: BE}/\sqrt{n}$, by Lemma \ref{l: Laplace transform} we
have
\begin{align}
    \P \left ( \norm{Ax} \le \D n \right )
    &\le (\Cref{l: Laplace transform}\D \sqrt{n})^n \notag
\intertext{and so,}
    \P \left ( \exists x \in \mathcal{M} \mid \norm{Ax} \le \D n
       \right )
    &\le |\mathcal{M}| (\Cref{l: Laplace transform}\D \sqrt{n})^n.
\label{norm on net}
\end{align}
We shall show that $Q_1$ can be chosen so that the last quantity
will be less than $(2e)^{-n}$. Recall that by Lemma \ref{l: net
singular}, there exists a $\D$-net $\mathcal{N}$ for the set of
vectors of $(\D,Q_1)$-singular profile satisfying
\[
  |\mathcal{N}|
  \le \left ( \frac{\Cref{l: net singular}}
                   {\D \sqrt{n}}
               \eta^{\cref{l: net singular}} \right )^n,
\]
provided
\begin{equation} \label{Q_1}
    C(\eta) <\Cref{d: regular profile} Q_1.
\end{equation}
Covering each cube of size $\D$ with the center in $\mathcal{N}$
by the cubes of size $\frac{\D}{2 \Cref{l: norm of A}}$ and using
Lemma \ref{l: volumetric}, we obtain
\[
  |\mathcal{M}|
  \le |\mathcal{N}| \cdot
     N(\D B_{\infty}^n, \D B_{\infty}^n,
           \frac{1}{2 \Cref{l: norm of A}})
  \le \left ( \frac{6\Cref{l: net singular} \cdot
                      \Cref{l: norm of A}}
                   {\D \sqrt{n}}
               \eta^{\cref{l: net singular}} \right )^n.
\]
Substitution of this estimate into \eqref{norm on net} yields
\begin{align*}
   \P \left ( \exists x \in \mathcal{N} \mid \norm{Ax} \le \D n
       \right )
  &\le \left ( \frac{6\Cref{l: net singular} \cdot
                      \Cref{l: norm of A}}
                   {\D \sqrt{n}}
               \eta^{\cref{l: net singular}} \right )^n
  \cdot (\Cref{l: Laplace transform}\D \sqrt{n})^n \\
  &\le (C' \eta^{\cref{l: net singular}})^n.
\end{align*}
Now choose $\eta$ so that $C' \eta^{\cref{l: net singular}}< 1/e$
and choose $Q_1$ satisfying \eqref{Q_1}. With this choice the
probability above is smaller than $e^{-n}$. Combining this
estimate with Lemma \ref{l: norm of A}, we have that $\norm{A} \le
\Cref{l: norm of A} \sqrt{n}$ and $\norm{Ax} \ge \D n$ for all $x
\in \mathcal{N}$ with probability at least $1-2e^{-n}$.

Let $y \in V_S$ be a vector of $(\D,Q_1)$-singular profile. Choose
$x \in \mathcal{N}$ such that $\norm{x-y}_{\infty} \le
\frac{\D}{2\Cref{l: norm of A} }$. Then $\norm{x-y} \le \frac{\D
\sqrt{n}}{2\Cref{l: norm of A} }$ and
\[
  \norm{Ay} \ge \norm{Ax} - \norm{A(x-y)}
  \ge \D n - \norm{A}\norm{x-y}
  \ge \frac{\D}{2} n.
\]

{\bf Case 2} Assume that $\Cref{t: singular profile} n^{-3/2} \le
\D < \D_1 = \cref{l: BE}/n$. Let $\Om_1$ be the event that
$\norm{Ax} < \frac{\D_1}{2} n = \cref{l: BE} / 2$ for some vector
of $(\D_1,Q_1)$-singular profile. We proved in Case 1 that
\begin{equation}  \label{Om_1}
  \P(\Om_1) < 2e^{-n}.
\end{equation}
Let $Q_2>1$ be a number to be chosen later and let $W$ be the set
of all vectors of $(\D_1, Q_1)$-singular and $(\D, Q_2)$-regular
profile. By Lemma \ref{l: regular disrepancy} any  vector $x \in
W$ satisfies
\[
  \P \left ( \left | \sum_{j=1}^n \b_j x_j \right | \le  t
     \right )
  \le  \Cref{l: regular disrepancy} Q_1 t
\]
for all $t \ge \D_1$.

Now we can finish the proof as in Case 1. Since $\D \sqrt{n} \ge
\D_1$, Lemma \ref{l: Laplace transform} implies
\[
    \P \left ( \norm{Ax} \le \D n \right )
    \le (C'\D \sqrt{n})^n
\]
for any $x \in W$. Here $C'=\Cref{l: Laplace transform}\cdot
\Cref{l: regular disrepancy} Q_1$.

Let $\mathcal{N}$ be the smallest $\frac{\D}{2 \Cref{l: norm of A}
}$-net in $W$ in $\ell_{\infty}$ metric. Note that $\D \ge
\Cref{t: singular profile} n^{-3/2} \ge \overline{\Cref{l: net
singular}} n^{-3/2}$. Arguing as in the Case 1, we show that
\[
  |\mathcal{N}|
  \le \left ( \frac{6\Cref{l: net singular}
                    \cdot \Cref{l: norm of A}}{\D \sqrt{n}}
       \eta^{\cref{l: net singular}} \right )^n
\]
for any $\eta$ satisfying
\begin{equation} \label{Q_2}
    C(\eta) <\Cref{d: regular profile} Q_2.
\end{equation}
 Hence,
\[
   \P \left ( \exists x \in \mathcal{N} \mid \norm{Ax} \le \D n
       \right )
    \le |\mathcal{N}| (C'\D \sqrt{n})^n
    \le (C''  \eta^{\cref{l: net singular}} )^n.
\]
Choose $\eta$ so that the last quantity is less than  $e^{-n}$ and
choose $Q_2$ so that \eqref{Q_2} holds. Then the approximation
argument used in Case 1 shows that the inequality
\[
  \norm{Ay}
  \ge \frac{\D}{2} n
\]
holds for any $y \in W$ with probability greater than $1-e^{-n}$.
Combining it with \eqref{Om_1}, we complete the prof of Case 2.
Finally, we unite two cases setting $Q_0=\max(Q_1,Q_2)$.
\end{proof}

\section{Proof of Theorem \ref{t: main}.}   \label{s: main proof}

To prove Theorem \ref{t: main} we combine the probability
estimates of the previous sections. Let $\e > \cref{t:
main}/\sqrt{n}$, where the constant $\cref{t: main}$ will be
chosen later. Define the exceptional sets:
\begin{align*}
  \Om_0 &= \{ \om \mid \norm{A} > \Cref{l: norm of A} \sqrt{n} \},
   \\
  \Om_P &= \{ \om \mid \exists \, x \in V_P \
                      \norm{Ax} < \Cref{l: peaked} \sqrt{n} \
                \}.
\end{align*}
Then Lemma \ref{l: norm of A} and Lemma \ref{l: all peaked} imply
\[
  \P(\Om_0)+\P(\Om_P) \le 3  \exp( - \cref{l: peaked} n ).
\]

Let $Q_0$ be the number defined in Theorem \ref{t: singular
profile}. Set
\[
  \D =\frac{\e}{2 \Cref{l: regular disrepancy} Q_0 \cdot n}.
\]
 The
assumption on $\e$ implies $\D \ge \Cref{t: singular profile}
n^{-3/2}$ if we set $\cref{t: main}=2 \Cref{l: regular disrepancy}
Q_0 \cdot \Cref{t: singular profile}$. Denote by $W_S$ the set of
vectors of $(\D, Q_0)$-singular profile and by $W_R$ the set of
vectors of $(\D, Q_0)$-regular profile. Set
\begin{align*}
  \Om_S &= \{ \om \mid \exists \, x \in W_S \
                      \norm{Ax} < \frac{\D}{2}n
                          =\frac{1}{4 \Cref{l: regular disrepancy} Q_0}
                             \e
                \}, \\
  \Om_R &= \{ \om \mid \exists \, x \in W_R \
                      \norm{Ax} <  \frac{\D}{2 \sqrt{n}}
                      =\frac{1}{4 \Cref{l: regular disrepancy} Q_0}
                        \e \cdot n^{-3/2}
                \}.
\end{align*}
By Theorem \ref{t: singular profile}, $\P(\Om_S) \le 3e^{-n}$, and
by Theorem \ref{t: regular profile}, $\P(\Om_R) \le \e/2$. Since
$S^{n-1} = V_P \cup W_S \cup W_R$, we conclude that
\[
  \P(\om \mid \exists \, x \in S^{n-1} \
                      \norm{Ax}
                      <\frac{1}{2 \Cref{l: regular disrepancy} Q_0}
                        \e \cdot n^{-3/2}
                \}
  \le \e/2 + 4 \exp( - \cref{l: peaked} n )
  <\e
\]
for large $n$. \qed
\begin{remark}
The proof shows that the set of vectors of a regular profile is
critical. On the other sets the norm of $Ax$ is much greater with
probability exponentially close to $1$.
\end{remark}


\begin{thebibliography}{GGM}
%
\bibitem[BY]{BY}
Z. D. Bai, Y. Q. Yin,{\em Limit of the smallest eigenvalue of a
large-dimensional sample covariance matrix},  Ann. Probab.  21
(1993),  no. 3, 1275--1294.

\bibitem[DS]{DS}
K. Davidson, S. J. Szarek, {\em Local operator theory, random
matrices and Banach spaces},  Handbook of the geometry of Banach
spaces, Vol. I,  317--366, North-Holland, Amsterdam, 2001.

\bibitem[E]{Es}
C. G. Ess\'{e}en, {\em On the concentration function of a sum of
independent random variables},  Z. Wahrscheinlichkeitstheorie und
Verw. Gebiete  9  1968 290--308.

\bibitem[Ha1]{Ha1}
G. Hal\'{a}sz, {\em On the distribution of additive arithmetic
functons}, Acta Arithmetica, XXVII (1975), 143--152.

\bibitem[Ha2]{Ha2}
G. Hal\'{a}sz, {\em Estimates for the concentration function of
combinatorial number theory and probability}, Per. Math. Hung.
8(3-4), (1977), 197--211.

\bibitem[KKS]{KKS}
J. Kahn, J. Koml\'{o}s, E. Szemer\'{e}di, {\em On the probability
that a random $±1$-matrix is singular}, J. Amer. Math. Soc.  8
(1995), no. 1, 223--240.

\bibitem[K1]{K1}
J. Koml\'{o}s, {\em On the determinant of $(0,\,1)$ matrices},
Studia Sci. Math. Hungar.  2  (1967), 7--21.

\bibitem[K2]{K2}
Koml\'{o}s, J. {\em On the determinant of random matrices},
Studia Sci. Math. Hungar.  3  (1968), 387--399.

\bibitem[KSC]{KSC}
 V. F. Kolchin, B. A. Sevast'yanov and V. P. Chistyakov,  {\em Random
allocations},  Scripta Series in Mathematics. V. H. Winston \&
Sons, Washington, D.C.; distributed by Halsted Press [John Wiley
\& Sons], New York-Toronto, Ont.-London, 1978.

\bibitem[LPRT]{LPRT}  A.~E.~Litvak, A.~Pajor, M.~Rudelson and
N.~Tomczak-Jaegermann, {\em Smallest singular value of random
matrices and geometry of random polytopes}, Adv. Math., to appear.
%
\bibitem[LPRTV1]{LPRTV1} A.~E.~Litvak, A.~Pajor, M.~Rudelson,
N.~Tomczak-Jaegermann and  R.~Vershynin, {\em Random Euclidean
embeddings in spaces of bounded volume ratio}, C. R. Acad. Sci.
Paris, S\'er. I Math.,
 339 (2004), 33--38.

 \bibitem[LPRTV2]{LPRTV2} A.~E.~Litvak, A.~Pajor, M.~Rudelson,
N.~Tomczak-Jaegermann and  R.~Vershynin, {\em Euclidean embeddings
in spaces of finite volume ratio via random matrices},  to appear
in J. Reine Angew. Math.

\bibitem[P]{P}
G. Pisier, {\em The volume of convex bodies and Banach space
geometry}, Cambridge Tracts in Mathematics, 94. Cambridge
University Press, Cambridge, 1989.

\bibitem[St]{Strook} D. Stroock, {\em Probability theory. An analytic
view}, Cambridge Univ. Press 1993.

\bibitem[Sz1]{Sz1}
S. J. Szarek, {\em Spaces with large distance to $l\sp n\sb
\infty$ and random matrices},  Amer. J. Math.  112  (1990),  no.
6, 899--942.

\bibitem[Sz2]{Sz2}
S. J. Szarek,{\em Condition numbers of random matrices}, J.
Complexity 7  (1991),  no. 2, 131--149.

\bibitem[TV1]{TV1}
T. Tao, V. Vu, {\em On random $\pm 1$ matrices: singularity and
determinant}, preprint.

\bibitem[TV2]{TV2}
T. Tao, V. Vu, {\em On the singularity probability of random
Bernoulli matrices}, preprint.


\end{thebibliography}
\end{document}